\documentclass[10pt]{amsart}
\usepackage{amssymb}
\usepackage[utf8]{inputenc}
\usepackage{array}

\newtheorem{theorem}{Theorem}[section]
\newtheorem{lemma}[theorem]{Lemma}
\newtheorem{corollary}[theorem]{Corollary}
\newtheorem{proposition}[theorem]{Proposition}
\newtheorem{remark0}[theorem]{Remark}
\newtheorem{example0}[theorem]{Example}
\newtheorem{definition}[theorem]{Definition}

\newenvironment{example}{\begin{example0}\rm}{\end{example0}}
\newenvironment{remark}{\begin{remark0}\rm}{\end{remark0}}

\newcommand{\defref}[1]{Definition~\ref{#1}}
\newcommand{\propref}[1]{Proposition~\ref{#1}}
\newcommand{\thmref}[1]{Theorem~\ref{#1}}
\newcommand{\lemref}[1]{Lemma~\ref{#1}}

\newcommand{\exref}[1]{Example~\ref{#1}}

\def\max{{\mathfrak{m}}}                   

\def\res{{\mathbf{k}}}

\def\HF{{\operatorname{H\!F}}}

\def\dim{\operatorname{dim}}
\def\ann{\operatorname{Ann}}

\def\Hom{\operatorname{Hom}}
\def\length{\operatorname{Length}}
\def\Max{\operatorname{Max}}

\begin{document}

\title[Almost finitely generated Inverse systems and reduced $\res-$algebras]{Almost finitely generated Inverse systems and reduced $\res-$algebras}

\author{J. Elias ${}^{*}$}
\thanks{${}^{*}$ Partially supported by  PID2022-137283NB-C22}
\address{Joan Elias
\newline \indent Departament de Matem\`{a}tiques i Inform\`{a}tica
\newline \indent Universitat de Barcelona (UB)
\newline \indent Gran Via 585, 08007
Barcelona, Spain}
\email{{\tt elias@ub.edu}}

\author{M. E. Rossi ${}^{**}$}
\thanks{${}^{**}$
Partially supported by
: PRIN-MIUR 2020355B8Y}

\address{Maria  Evelina  Rossi
\newline \indent Dipartimento di Matematica
\newline \indent Universit\`{a} di Genova
\newline \indent Via Dodecaneso 35, 16146 Genova, Italy}
\email{{\tt rossim@dima.unige.it}}


\begin{abstract}
The purpose of this paper is to characterize one-dimensional local domains, or more in general reduced, in terms of its Macaulay's inverse system. This leads to study almost finitely generated modules in the divided power ring.
We specialize the results to a numerical semigroup ring  by computing explicitly  its  inverse system.
In the graded case we characterize reduced arithmetically Gorenstein $0$-dimensional schemes.
Several examples are given.
\end{abstract}

\maketitle

\section{Introduction}

Let $(R, \max, \res)$ be a complete Noetherian local ring where $\res=R/\max$ denotes the residue field and let $E_R(\res)$ be the injective envelope.  Given an $R$-module $M, $  the  Matlis dual   $M^\vee=\Hom_R(M,E_R(\res))$
  defines a controvariant, additive and  exact functor from the category of the  $R$-modules into itself. In particular, by Matlis duality,   the functor  $(-)^\vee$ is  an anti-equivalence between the category of the finitely generated $R$-modules and the category of the Artinian $R$-modules, \cite[Theorem 3.2.13]{BH97}.

Let  $R$ be the power series ring $ \res[\![x_1, \dots, x_n]\!] $ or the polynomial ring
$ \res[x_1, \dots, x_n] $ over a field $\res.$ We denote by $\max$ the maximal  ideal of $R$ generated by $x_1,\dots, x_n. $ It is  known that the injective envelope $E_R(\res{})$ is isomorphic as $R$-module to the divided power ring $ \Gamma= \res_{DP} [y_1, \dots, y_n]  $ (see \cite{Gab59}, \cite{Nor72b}) also denoted in the literature by $ \res  [x_1^{-1}, \dots, x_n^{-1}]  $.

 Macaulay in \cite[IV]{Mac16} proved a particular case of Matlis duality, called Macaulay's correspondence,
between the
  ideals $I\subseteq R  $ such that $R/I$ is an Artinian local ring and $R$-submodules  $I^{\perp}$   of   $\Gamma $ of $R$ which are finitely generated. For more detailed information concerning Macaulay's Inverse System   see also \cite{EisGTM}, \cite{Ems78}, \cite{EI78}, \cite[Appendix A]{IK99}.    This theory  was   recently extended by Kleiman and Kleppe in \cite{KK22} to the more general situation of     $A$-algebras, where $A$ is any Noetherian ring.

 Macaulay's  correspondence is an effective method for computing Gorenstein Artinian rings, see  \cite{CI12}, Section 1, \cite{Iar94}, \cite{Ger96} and \cite{IK99}.
 An Artinian  Gorenstein $\res$-algebra $A=R/I$ of socle degree $s$   corresponds to a cyclic $R$-submodule   of $\Gamma$ generated by a polynomial $F\neq 0$ of degree $s$.

 The authors extended  Macaulay's correspondence to   $d$-dimensional local Gorenstein $\res$-algebras proving that they are in bijective  correspondence with  suitable submodules of $\Gamma$, called $G$-admissible, see \cite{ER17}. This  result was extended by S. Masuti, P. Schulze and L.Tozzo to any Cohen-Macaulay $\res$-algebra (local or graded), given   the dimension of the socle of a suitable Artinian reduction, see   \cite{MT18}, \cite{ST19}.

Despite several advances on the topic, many basic problems remain open. For instance, in the Artinian case,  it is not known which polynomials $F$  in $\Gamma$ correspond to the ideals of R which are a complete intersection, see \cite[pag. 261]{IK99}. In positive dimension, it would be interesting to describe Macaulay's inverse system of domains and, more in general, of reduced $\res$-algebras.
The last question is the   purpose of this paper.

\vskip 2mm
We   briefly describe the organization of the paper. In Section 2 we present the   main tools concerning Macaulay's Inverse System that will be used in the paper. In Section 3,  inspired by results of W.D. Weakley \cite{Wea83},  we prove that the inverse system of a  $1$-dimensional local domain    is an almost finitely generated  $R$-submodule of $\Gamma  $ and conversely, see Theorem \ref{prime}.  As a consequence, taking advantage of the properties of the $G$-admissible $R$-submodules of $\Gamma,$ see Definition \ref{G-modules}, we characterize the inverse system of  one dimensional local Gorenstein domains, see Proposition \ref{localdomains}. In Section 4, \thmref{semigroup},  we   explicitly describe the generators of the almost finitely generated dual module   of a numerical semigroup ring,  see also \cite[Theorem 2.1]{EW21} for analogous results recently  obtained  by  K. Eto and K. Watanabe. The defining ideal of a numerical semigroup ring is an homogeneous ideal in a  weighted polynomial ring, hence in the process we need to study the inverse system in the non-standard graded case.
In Section 5 we study the inverse system of one dimensional reduced $\res$-algebras, see Proposition \ref{radical}.  In the graded setting, this is the case of the coordinate ring of a set $X$ of distinct points in the projective space. In Theorem \ref{ger-link} we present conditions on the inverse system of any Gorenstein zero-dimensional scheme $X$  for being a reduced scheme.   This results completes Theorem 3.14 in \cite{ER22}.
In particular we translate  the problem   in terms of the identifiability (in the sense of \cite{AC20}) of a specific polynomial in the inverse system of $X.$

 \vskip2mm
Throughout the paper several examples are given. The computations are performed  by using the computer algebra system  Singular  \cite{DGPS}  and in particular the Singular library INVERSE-SYST.lib, \cite{E-InvSyst14}.


\section{Preliminaries}

We recall   that $\Gamma =\res_{DP} [y_1, \dots, y_n]$ is an  $R$-module acting  $R$ by    {\it {contraction}} denoted   by $\circ$.
Given  $\alpha, \beta\in \mathbb N^n$, we  denote by $ x^{\alpha}= x_1^{\alpha_1} \cdots x_n^{\alpha_n} \in R$ and $ y^{\beta}= y_1^{\beta_1} \cdots y_n^{\beta_n}\in \Gamma$,  then
$$
x^{\alpha} \circ y^{\beta} =
\left\{
\begin{array}{ll}
 y^{ \beta-  \alpha}   & \text{ if }  \beta_i \ge \alpha_i \text{ for all  } i=1,\dots, n\\ \\
0 & \text{ otherwise}
\end{array}
\right.
$$
If the characteristic of field $\res=R/\max$ is zero then there is a natural isomorphism of $R$-algebras between $(\Gamma,\circ)$ equipped with an internal product and the polynomial ring  replacing the contraction with the partial derivatives.
This action is sometimes called the “apolarity” action of $R$ on $\Gamma$
defined by
$$
\begin{array}{ cccc}
\circ: & R  \times \Gamma  &\longrightarrow &  \Gamma   \\
                       &       (f , g) & \to  &  f  \circ g = f (\partial_{y_1}, \dots, \partial_{y_n})(g)
\end{array}
$$

\noindent
where $  \partial_{y_i} $ denotes the partial derivative with respect to $y_i.$

In this paper we assume $\res$ of any characteristic, unless otherwise specified.

In both cases (contraction or derivation) the action of $R$ on $\Gamma$  lowers degree. Thus, $\Gamma$ is not a finitely generated
R-module.     Let $R_i$  (resp.    $\Gamma_i$) be the sub-$\res$-vector space of $R$ (resp. of  $\Gamma$)  generated by the standard monomials of degree $i\ge 0.$   Notice that the apolarity action induces a non-singular $\res$-bilinear pairing:
$$
\begin{array}{ cccc}
\circ: & R_j  \times \Gamma_j   &\longrightarrow &  \res
\end{array}$$
for every $j \ge 0.$ The standard grading  will be extended to a weighted grading  in Section 3. However, unless explicitly stated otherwise, a polynomial ring has the standard grading in which all the variables have weight $1.$

 \vskip 2mm
 If $I\subset R$ is an ideal of $R, $ then  $(R/I)^\vee={\rm Hom}_R(R/I,E_R(\res))$ is the $R$-submodule of $\Gamma$
$$
{I^{\perp}} =\{ F \in \Gamma \ |\  I \circ F = 0 \  \}.
$$
This submodule of $\Gamma $  is called the {\it{Macaulay's inverse system of $I$.}} Remark that $I^{\perp}$ is also an $R/I$-module since $I \circ I^{\perp} =0$.

\noindent Conversely, given an $R$-submodule $W$ of $\Gamma, $  the dual $W^\vee={\rm Hom}_R(W,E_R(\res))$ is the ring $R/\ann_R(W)$ where
$$ \ann_R(W)= \{ g \in R \ \mid \ g \circ F= 0 \ \mbox{ for \ all \ } g \in W\}
$$
 is an  ideal of $R$. Macaulay's correspondence  in \cite[IV]{Mac16} gives a correspondence
between the
  ideals $I\subseteq R  $ such that $R/I$ is an Artinian local ring and $R$-submodules     of   $\Gamma $ of $R$ which are finitely generated. In particular   Macaulay proved that   Artinian  Gorenstein $\res$-algebras $A=R/I$ of socle degree $s$   correspond to cyclic $R$-submodules  of $\Gamma$ generated by a polynomial $F\neq 0$ of degree $s$.


 \vskip 2mm
In \cite{ER17} the authors extended Macaulay's correspondence to the $d$-dimensional Gorenstein  $\res$-algebras $R/I$ characterizing the generators of
the $R$-submodules of $\Gamma$ (not finitely generated), called $G$-admissible, in correspondence with $R/I. $ We present here the definition of $G$-admissible for $d=1$.

\begin{definition}
\label{G-modules}  An   $R$-submodule $M $ of  $\  \Gamma  $    is called {\bf{G-admissible}} if it admits a countable system of  generators $\{H_l\}_{ l\in \mathbb N_+} $ satisfying  the following conditions
\begin{enumerate}
\item[(1)] There exists a linear form   $z \in R$ such that for all $l\in \mathbb N_+$
$$
z \circ H_l =
\left\{
\begin{array}{ll}
 H_{l-1}   & \text{ if }  l> 1 \ \\
0 & \text{ otherwise.}
\end{array}
\right.
$$
\item[(2)] $\ann_R(\langle H_l\rangle)\circ H_{l+1}=\langle H_{1}\rangle$ for all $l\in \mathbb N_+$.
\end{enumerate}

\noindent
If this is the case, we say that $M = \langle H_l, l\in \mathbb N_+\rangle$ is a
$G$-admissible $R$-submodule of $\Gamma$ with respect to the linear form   $z \in R$.
\end{definition}


We present the main result of \cite{ER17} in the one-dimensional case.

\medskip
\begin{theorem}[\cite{ER17}, Theorem 3.8]
\label{bijecGor}
There is a one-to-one correspondence $\mathcal C$ between the following sets:

\noindent
(i)
 one-dimensional  Gorenstein $\res$-algebras $A=R/I$,

\noindent
(ii) non-zero $G$-admissible $R$-submodules $M=\langle H_l, l\in \mathbb N_+\rangle$ of $\Gamma$.

\noindent
In particular, given an ideal $I\subset R$ with $A= R/I$ satisfying $(i)$ and $z$ a linear regular element  modulo $I, $   then
 $$\mathcal C(A)= I^{\perp} = \langle H_l,  l\in \mathbb N_+\rangle \subset S   \ \ {\text{with}}\ \
  \langle H_l\rangle=(I+(z^l))^\perp $$ is $G$-admissible.
 Conversely, given an  $R$-submodule $M$ of $\Gamma $ satisfying  (ii), then $$\mathcal C^{-1}(M)=R/I  \ \ {\text{with}}\ \
 I= \ann_R(M) = \bigcap_{l\in \mathbb N_+}\ann_R(\langle H_l\rangle).$$
\end{theorem}

\bigskip
\section{Inverse systems of local domains}

We recall the definition of almost finitely generated $A$-module, where $A$ is a  Noetherian  local ring, see \cite{Wea83}, \cite{HL85}.

\begin{definition}
An $A$-module $M$ is {\bf{almost finitely generated}}, a.f.g. for short, if $M$ is not finitely generated, but
any proper sub-$A$-module $N$ of $M$ is finitely generated.
Moreover, $M$ is a {\bf{divisible}} $A$-module if for all  $a\in R\setminus \{0\}$  it holds $M=a M$.
\end{definition}

As before, let  $R$ be   the power series ring $R= \res[\![x_1, \dots, x_n]\!] $.
The following remark is an straightforward consequence of the definition of a.f.g. $R$-module that will be used systematically along this section.

\medskip
\begin{lemma} \label{lemmaprep}  Let $I\varsubsetneq J $ be  ideals of $R$ such that  $\dim R/J >0$. Then $I^{\perp}$ is not an a.f.g. $R$-module.
\end{lemma}
\begin{proof}
Notice that $J^{\perp}$ is not finitely generated because $R/J$ is not Artinian.
Since $J^{\perp}\varsubsetneq I^{\perp}$, then $I^{\perp}$ is not an a.f.g. $R$-module.
\end{proof}

\medskip
 In the next result we characterize prime ideals in terms of their inverse systems. We remark that $I^{\perp}$ is an $R$-module, but also an $R/I$-module since $I \subseteq Ann_R(I^{\perp}). $
The following  result can eventually be  deduced from \cite[Proposition 2.6]{Wea83}, here we give a short and straight proof.

\medskip
\begin{theorem}
\label{prime}
Let $I$ be an ideal of $R$   such that  $A=R/I$ is one-dimensional.
The following conditions are equivalent:
\begin{enumerate}
  \item[(i)] $I$ is a prime ideal,
  \item[(ii)] $I^{\perp}$ is an a.f.g. $A$-module,
  \item[(iii)] $I^{\perp}$ is a divisible $A$-module.
\end{enumerate}
\end{theorem}
\begin{proof}
Assume that $I$ is a prime ideal.
Since $R/I$ is not an Artinian ring,  $I^{\perp}$ is not a finitely generated $R$-module.
Let $N\subsetneq I^{\perp}$ be an $R$-submodule of $I^{\perp}$.
Since $I$ is a prime ideal and  $I\subsetneq  J:=\ann_R(N)$, we have that $R/J$ is Artinian and hence
  $J^{\perp}=N$ is finitely generated. Hence $I^{\perp} $ is an a.f.g. $R$-module and also a.f.g. $A$-module.

Assume that $I^{\perp}$ is an a.f.g. $A$-module. From \cite{Wea83}, Proposition 1.1 (3), we get that $I^{\perp}$ is divisible.

If $I^{\perp}$ is a divisible $A$-module then for all $a\in A\setminus \{0\}$ it holds { {$I^{\perp}=a\circ I^{\perp}$}}, i.e. the morphism of $R$-modules
$$
I^{\perp} \overset{a\circ}{\longrightarrow} I^{\perp}
$$
is an epimorphism.
By Matlis duality we deduce that the product by $A \overset{a.}{\longrightarrow} A$ is a monomorphism, i.e. $a$ is a non-zero divisor.
\end{proof}

\medskip
\begin{remark}
Notice that if $I^{\perp}$ is a divisible $A$-module, then $I^{\perp}=\max\circ I^{\perp}$.
This condition does not imply $I^{\perp}=0$ since
$I^{\perp}$ is  not, in general, a finitely generated $A$-module, hence Nakayama Lemma does not hold in this situation.
\end{remark}

\medskip
In the next result we characterize the $G$-admissible sets  which are divisible.

\medskip
\begin{lemma}
\label{gens-inv}
Let $A=R/I$ be a one-dimensional Gorenstein ring.
Let $\{H_t; t\ge 1\}$ be a G-admissible system of generators of $I^{\perp}$.
Then for all $H\in I^{\perp}$ there exists $r\in \mathbb N$ and $F\in R$ such that
$H=F \circ H_r$.
\end{lemma}
\begin{proof}
It follows easily from  the condition $(1)$ of G-admissibility, \defref{G-modules}, because if $t > r, $ then $H_t= z^{t-r} \circ H_r$.
\end{proof}

\medskip
\begin{proposition}\label{localdomains}
Let $I$ be an ideal of $R$   such that  $A=R/I$ is one-dimensional and Gorenstein. Let
$I^{\perp}=\langle H_t; t \ge 1\rangle$, where $\{H_t; t\ge 1\}$  is  a G-admissible system of generators with respect to $z\in R$.
 Then $I$ is prime if and only if  for all $a\in { {A}}\setminus \{0\}$      and for all $t\ge 1$
there exist  $F\in R$ and $r\ge t$ such that
$$
H_t= a\circ (F\circ H_r).
$$
\end{proposition}
\begin{proof}
Assume that $I$ is a prime ideal.
From Theorem \ref{prime},  $I^{\perp}$ is a divisible $A$-module.
Then for all $a\in A\setminus \{0\}$ and for all $t\ge 1$ there exist $G\in I^{\perp}$ such that $H_t=a \circ G$.
From  \lemref{gens-inv} there exists $F\in R$ such that $G=F\circ H_r$ for some $r\in \mathbb N$.
Hence
$$
H_t= a\circ G= a\circ (F\circ H_r).
$$
It is easy to prove that $r\ge t$. Let $\{H_t; t\ge 1\}$ a  G-admissible system of generators  with respect to $z$.
If $r< t$ then
$$
H_1= z^{t-1}\circ H_t= (F a)\circ (z^{t-1}\circ H_r)= (F a)\circ 0=0
$$
which is not possible.

Assume now that for all $a\in A\setminus \{0\}$  and for all $t\ge 1$
there exists $F\in R$ and $r\ge t$ such that
$$
H_t= a\circ (F\circ H_r).
$$

From Theorem \ref{prime}  we only have to prove that $I^{\perp}$ is divisible, i.e.
for all $a\in A\setminus \{0\}$ and $H\in I^{\perp}$ there exist $L\in I^{\perp}$ such that $H=a\circ L$.

Let $0\neq a\in A$ be an element of $A$   and $H\in I^{\perp}$.
Then there exist $q\in R$ and $w\in \mathbb N$ such that
$$
H= q \circ H_{w}.
$$
From the hypothesis there exist $F\in R$ and $r\ge w$ such that
$H_{w}= a\circ (F\circ H_r)$,
so
$$
H=q\circ H_{w}= a\circ( (qF)\circ H_r)=a\circ L.
$$
\end{proof}

\begin{remark} Let $A=R/I$ be a local ring with maximal ideal $\max$.
We recall that the  valuation with respect to $\max$ is the function
$$
\begin{array}{cccl}
v_{\max}:& A\setminus \{0\}&\longrightarrow & \mathbb N\\
&a &\mapsto &
v_{\max}(a)=\Max\{u\mid a\in\max^u\}
\end{array}
$$
and $v_{\max}(0)=+\infty$.  It is easy to see that in Proposition  \ref{localdomains}, instead of all the elements $a\in A, $ we may restrict the control to the elements $a \in A $ such that  $v_{\max}(a)\le e-1$ where    $e$ is the  multiplicity of $A=R/I$.
In fact, since $A$ is a one-dimensional Cohen-Macaulay local ring
we have
$$
\max^{e-1+u}=z^u{\max^{e-1}}.
$$
for every $u \ge 1.$ Hence if  $v_{\max}(a)= e-1 + u > e-1, $ then $a=z^u c$ with $c\in \max^{e-1}\setminus\max^{e}$, in particular $v_{\max}(c)\le e-1$. In the proof it is enough to take
$L=F\circ H_{r+u}$ then we have
$$
a\circ L=a\circ (F\circ H_{r+u})= (z^u c)\circ (F\circ H_{r+u})
=c\circ (F\circ H_r)=H.
$$
\end{remark}

\medskip
\section{Inverse system of numerical semigroup rings}

Let $1\le a_1 \le \cdots \le  a_n$ be an $n$-ple of positive integers and let $\omega=(a_1, \dots, a_n)$. Consider    the corresponding  ring $R=\res[\![ x_1, \dots, x_n ]\!] $   where   $\deg_{\omega}(x_i)=a_i$, $i=1,\dots,n.$
Denote    $\omega(b_1,\dots,b_n):=\sum_{i=1}^n b_i a_i. $ If $x^K=x_1^{k_1}\dots x_n^{k_n}, $ then
$\deg_{\omega}(x^K)=\omega (K) $ for $K=(k_1, \dots, k_n) \in \mathbb N^n. $
 Denote by $R_{\omega, j}$  (resp. $\Gamma_{\omega, j}$)  the sub-$\res$-vector space of $R$ (resp. $\Gamma$) generated by the monomials    of degree $j$    with respect to $\omega$.

Then    the apolarity action induces a non-singular $\res$-bilinear pairing:
$$
\begin{array}{ cccc}
\circ: & R_{\omega, j}  \times \Gamma_{\omega, j}  &\longrightarrow &  \res
\end{array}$$
for every $j \ge 0$, where  $\circ$ is the contraction.
As for the standard graded  case we have that if $I$ is an homogeneous ideal in the $\omega$-weighted ring  $R, $
then $I^\perp$ is homogeneous in the $\omega$-weighted  divided power ring $\Gamma $ and
$I^\perp=\oplus_j I^{\perp}_{\omega, j}$
where $I^{\perp}_{\omega, j}=\{F\in \Gamma_{\omega, j} \mid g\circ F=0 \text{  for all } g\in I\}$,
i.e. $(I^\perp)_{\omega, j}=I^{\perp}_{\omega, j}$. In fact we can repeat the same proof as in the standard case,  see \cite{Ger96}, Proposition 2.5.

\bigskip
We consider now local rings defined by
numerical semigroup rings and we   compute their inverse system.    Notice a description of the inverse system of numerical semigroup rings was also the main task of     \cite[Theorem 2.1]{EW21}. Here we prove the result in an easier way and we observe, following the definition, that it is an almost finitely generated $R$-module.

\medskip

Given  the  integers $1\le a_1< \cdots < a_n$ with $\gcd(a_1,\dots, a_n)=1, $ we denote by
$A(a_1,\dots,a_n):=R/I(a_1,\dots,a_n)$ the ring associated to the monomial curve with parameterization
$$
\begin{array}{ cccc}
\phi: &  R &\longrightarrow &  \res[\![t]\!]  \\
                       &       x_i & \to  &  t^{a_i}
\end{array}
$$
i.e. $I(a_1,\dots,a_n)=\ker(\phi)$.
If we denote   $\omega= (a_1, \dots, a_n) $   as before,  then the ideal
$I(a_1,\dots,a_n)$  is generated by the binomials
$x^K-x^L$ with $K, L\in \mathbb N^n$ such that $\deg_{\omega}(K)=\deg_{\omega}(L)$, see for instance \cite[Lemma 4.1]{Stu96}.

We denote by $\mathcal J$ the set of $j\ge 0$ such that there exists
$K\in \mathbb N^n$ with $\deg_{\omega}(K)=j$.
For all $j\in \mathcal J$
we define the following homogeneous form with respect to $\omega$ as follows:
$$
L_{\omega, j}=\sum_{\deg_{\omega}(K)=j} y^K  .
$$
Notice that $x_i\circ L_{\omega, j}=L_{\omega, j-a_i}$ for all  $j\in \mathcal J$ and $i=1, \dots, n$.

\medskip
\begin{theorem}
\label{semigroup}
Given the  integers $1\le a_1< \cdots < a_n$ with $\gcd(a_1,\dots, a_n)=1, $ then
$$I(a_1,\dots,a_n)^{\perp}=\bigcup_{j\in \mathcal J} \langle L_{\omega,j}\rangle$$ and the $R$-module
$I(a_1,\dots,a_n)^{\perp}$ is an a.f.g. $R$-module.
\end{theorem}
\begin{proof}
For all $j\in \mathcal J$, by the non-singular $\res$-bilinear
pairing induced by $\circ$ we get
$$
\dim_{\res} I^{\perp}_{\omega, j}=\dim_{\res} A(a_1,\dots,a_n)_{\omega, j}{=\dim_{\res}(\res[\![t]\!]_{j})=1},
$$
so we only have to prove that $L_{\omega,j}\in I^{\perp}_{\omega, j}$.
This is equivalent to show that for every $j \in \mathcal J, $
$$
(x^K-x^L)\circ L_{\omega,j}=0
$$
for all $K, L\in \mathbb N^n$ such that $\deg(K)=\deg(L)$.
Furthermore, this is equivalent to prove that for all $y^\alpha\in \Gamma_{\omega,j}$ such that
$x^K\circ y^\alpha\neq 0$ there exists a unique $y^\beta \in \Gamma_{\omega,j}$ with
$x^K\circ  y^\alpha= x^L\circ y^{\beta}.$
Since $x^K\circ y^\alpha\neq 0$ we have $\alpha-K\in \mathbb N^n$, so $\beta=L+(\alpha-K)\in \mathbb N^n$.

Since $I(a_1,\dots,a_n)$ is a prime ideal of $R$, by \thmref{prime}
the $R$-module
$I(a_1,\dots,a_n)^{\perp}$ is an a.f.g. $R$-module.
\end{proof}

\medskip
\begin{corollary}
\label{IX}
Given the  integers $1\le a_1< \cdots < a_n$ with $\gcd(a_1,\dots, a_n)=1, $ it holds
$$
(I(a_1,\dots,a_n)+(x_1^t))^{\perp}=\langle L_{\omega,j}\mid \deg_{y_1}(L_{\omega,j})\le t-1\rangle
$$
for all $t\ge 1$.
\end{corollary}
\begin{proof}
By Matlis' duality we know that
$$
(I(a_1,\dots,a_n)+(x_1^t))^{\perp}=I(a_1,\dots,a_n)^{\perp}\cap (x_1^t)^{\perp}.
$$
Being $(x_1^t)$ monomial, $(x_1^t)^{\perp}$ is homogeneous with respect the grading defined by $\omega$.
Hence
$$
(I(a_1,\dots,a_n)+(x_1^t))^{\perp}_{\omega ,j}=I(a_1,\dots,a_n)^{\perp}_{\omega ,j}\cap (x_1^t)^{\perp}.
$$
From this identity and the previous result we get the claim.
\end{proof}

\medskip
\begin{example}
Let us consider the integers $a_1=5, a_2=6, a_3=9$.
The corresponding monomial curve $C$ is  a complete intersection, so $C$ is  Gorenstein.
Moreover, the ideal defining  $C$ is minimally generated by two elements: $I(5,6,9)=(x_1^3-x_2x_3,x_2^3-x_3^2)$.
Then by using Corollary \ref{IX} we get:

\medskip
\begin{center}
\begin{tabular}{|c|l|c|}
  \hline
$I=I(5,6,9)$ & generator $L_{\omega,j}$ & $j$   \\ \hline\hline
$(I+(x_1))^{\perp}$ & $y_3^2+y_2^3$& $j=18$\\ \hline
$(I+(x_1^2))^{\perp}$ & $y_1y_3^2+y_1y_2^3$& $j=23$\\ \hline
$(I+(x_1^3))^{\perp}$ & $y_1^2y_3^2+y_1^2y_2^3$& $j=28$\\ \hline
$(I+(x_1^4))^{\perp}$ & $y_2y_3^3+y_2^4y_3+y_1^3y_3^2+y_1^3y_2^3$& $j=33$\\
  \hline
\end{tabular}
\end{center}
\end{example}

\medskip
\begin{example}
\label{non-gor}
Let us consider the integers $a_1=5, a_2=6, a_3=7$.
The corresponding monomial curve $C$ is not a complete intersection and  the Cohen-Macaulay type of $C$ is two.
The ideal defining  $C$ is minimally generated by three elements: $I(5,6,7)=(x_1^4-x_2x_3^2,x_2^2-x_1x_3, x_1^3x_2-x_3^3)$.
Then by using Corollary \ref{IX} we get:

\medskip
\begin{center}
\begin{tabular}{|c|l|c|}
  \hline
$I=I(5,6,7)$ & generators $L_{\omega,j}$ & $j$  \\ \hline\hline
$(I+(x_1))^{\perp}$ & $y_2y_3, y_3^2$& $j=13, 14$\\ \hline
$(I+(x_1^2))^{\perp}$ & $y_2^3+y_1y_2y_3, y_2^2y_3+ y_1y_3^2$& $j=18, 19$\\ \hline
$(I+(x_1^3))^{\perp}$ & $y_1y_2^3+y_1^2y_2y_3, y_2^4+y_1^2y_2y_3+ y_1^2y_3^2$& $j=23, 24$\\ \hline
$(I+(x_1^4))^{\perp}$ & $y_3^4+y_1^2y_2^3+y_1^3y_2y_3, y_1y_2^4+y_1^2y_2^2y_3+ y_1^3y_3^2$& $j=28, 29$\\
  \hline
\end{tabular}
\end{center}
\end{example}

\medskip
\begin{example}
Let us consider the integers $a_1=6, a_2=7, a_3=11, a_4=15$.
The corresponding monomial curve $C$ is not a complete intersection: the ideal defining  $C$ is minimally generated by five elements: $I(6,7,11,15)=(x_4^2-x_1^2x_2x_3, x_3x_4-x_1^2x_2^2, x_1x_4-x_2^3, x_3^2-x_2x_4, x_2x_3-x_1^3)$.
Then:
\medskip
\begin{center}
\begin{tabular}{|c|m{7cm}|c|}
  \hline
$I=I(6,7,11,15)$ & generators $L_{\omega,j}$ & $j$  \\ \hline\hline
$(I+(x_1))^{\perp}$ & $y_2^2, y_3^2+y_2y_4$& $j=14, 22$\\ \hline
$(I+(x_1^2))^{\perp}$ & $y_1y_2^2, y_1y_3^2+y_1y_2y_4+y_2^4$& $j=20, 28$\\ \hline
$(I+(x_1^3))^{\perp}$ & $y_3y_4+y_1^2y_2^2, y_1^2y_3^2+y_1^2y_2y_4+y_1y_2^4$& $j=26, 34$\\ \hline
$(I+(x_1^4))^{\perp}$ & $y_1y_3y_4+y_2^3y_3+y_1^3y_2^2, y_2y_3^3+y_2^2y_3y_4+y_1^3y_3^2+y_1^3y_2y_4+y_1^2y_2^4$& $j=32, 40$\\
  \hline
\end{tabular}
\end{center}

\noindent
Hence $A(6,7,11,15)/(x_1)$ is level of Cohen-Macaulay type two.
In \cite{MT18}, Example 1, it is shown that $A(6,7,11,15)/(x_1)$ is level and $A(6,7,11,15)/(x_1+x_2)$ does  not.
Notice that $x_1$ is homogeneous with respect to $\omega$ but $x_1+x_2$ is not homogeneous with respect to $\omega$.
\end{example}

\medskip
\section{Inverse systems of reduced rings}

As a consequence of Theorem \ref{prime}, we present  a   characterization of  the radical ideals $I \subseteq R=\res[\![x_1,\dots,x_n]\!]$
in terms of a.f.g. modules in the case $\dim R/I =1$.

\medskip
\begin{proposition}
\label{radical}
Let $I$ be an ideal of $R$ such that $\dim R/I=1$.
The following conditions are equivalent:
\begin{enumerate}
  \item[(i)] $I$ is a radical ideal,
  \item[(ii)] there  exist  $M_1,\dots, M_r$ a.f.g. sub-$R$-modules of $I^{\perp}$   such that
  $I^{\perp}=M_1+\cdots+M_r$.
\end{enumerate}
\end{proposition}
\begin{proof}
Assume that $I$ is a radical ideal.
Then $I=p_1\cap\cdots \cap p_r$ where $p_1,\dots, p_r$ are   prime ideals of $R$ such that $\dim R/p_i =1$ for all $i.$
Then the modules $M_i=p_i^{\perp}$ are a.f.g by \thmref{prime}.
and  $I^{\perp}=   M_1+\cdots+M_r. $  Hence  (ii) is proved.

Assume now (ii).
We define $J_i=\ann(M_i)$; by \thmref{prime} $J_i$ is a prime ideal.
Since $I\subset J_i$ the ideal $J_i$ has height at least $n-1$.
On the other hand, $M_i$ is not f.g. so $J_i$ is a height $n-1$ prime ideal of $R$.
By Matlis duality we get
$I=J_1\cap\cdots \cap J_r$, so $I$ is radical.
\end{proof}

In the following we give an example of a reduced   1-dimensional local ring obtained as a monomial curve linked to a straight line.

\begin{example}
Let us consider the non-Gorenstein monomial curve $C$ of \exref{non-gor};
its defining  ideal is $I_1=(x_1^4-x_2x_3^2,x_2^2-x_1x_3, x_1^3x_2-x_3^3)\subset R=\res[\![x_1,x_2,x_3]\!]$.
Let $L$ be the straight line defined by the ideal $I_2=(x_2,x_3)\subset R$.

The union $D=C\cup L$ is a complete intersection defined by the radical ideal
$J=I_1\cap I_2=(x_3^3-x_1^3x_2, x_2^2-x_1x_3)$; in other words: $C$ is linked to $L$ by the complete intersection $D$.

We know that
$J^{\perp}=I_1^{\perp}+I_2^{\perp}$.
Since $I_1, I_2$ are prime ideals, from \thmref{prime} the $R$-modules $I_1^{\perp}$
and $I_2^{\perp}$ are a.f.g., see \propref{radical}.

Notice that   $I_1^{\perp}$ is generated by $L_{\omega, j}$, $j\ge 0$, with $\omega=(5,6,7)$, see Theorem \ref{semigroup} and \exref{non-gor}.
By a straightforward computation we get that  $I_2^{\perp}$ is generated by $y_1^n$, $n\ge 0$.
The ring  $R/J$ is Gorenstein, so by using Singular library \cite{E-InvSyst14}, $J^{\perp}$ can be  generated by the following  $G$-admissible system of generators with respect to $x_1$:
$H_1=y_2y_3^2,  H_2=y_1 H_1+ y_2^3 y_3, H_3=y_1 H_2+ y_2^5, H_4=y_1 H_3+ y_3^5,
\dots$.

Next, it could be interesting to see explicitly how some generators of $J^{\perp}$ belong to $I_1^{\perp}+I_2^{\perp}$ and vice-versa.
For instance, let us consider the generator $H_2$ of $J^{\perp}$:
$$
H_2=
x_3^2\circ(y_2^3x_3^3+x_1x_2x_3^4+y_1^3y_2^4+y_1^4y_2^2y_3+y_1^5y_3^2)-y_1^5
=x_3^2\circ L_{\omega, 39}-y_1^5
\in I_1^{\perp}+I_2^{\perp}.
$$
\noindent
On the other hand let us consider the generator $L_{\omega,23}$ of $I_1^{\perp}$:
$$
L_{\omega,23}=y_1y_2^3+y_1^2y_2y_3=x_3\circ H_3\in J^{\perp},
$$
and the generator $y_1^3$ of $I_2^{\perp}$:
$$
y_1^3=x_2x_3^2\circ H_4\in J^{\perp}.
$$
\end{example}

\medskip
\bigskip
From now on, $R$ will be  the polynomial ring   $  \res[x_1, \dots, x_n]$ with the standard grading (deg $ x_i=1$).
In this part we use results in \cite{Ger96} and \cite{ER22}, hence we from now on we assume  $\res$ is an algebraically closed field of characteristic zero and {\bf{$\circ$ is the derivation}}.

\medskip
Let $X$  be a zero-dimensional scheme, the first purpose is to understand when  $X$ is reduced, that is $X$ consists of   a set of distinct points.
Recall that $I(X)^{\perp} $ is  not a finitely generated graded $R$-module   and $(I(X)^{\perp})_j= (I(X)_j)^{\perp} $ for every $j \ge 0.$
Moreover the Hilbert function of $R/I(X)$ is
$$\HF_{R/I(X)}(j)= \dim_{\res}(R_j/I(X)_j)= \dim_{\res} (I(X)^{\perp})_j,$$
$j\ge 0$, see \cite[Proposition 2.5]{Ger96}.

Let $z$   be a linear form in $R$ such that $z(P_i) \neq 0$ for every $i=1,\dots, r,$ that is $z$ is a not zero divisor in $R/I(X).$ Then we will say that $R/I(X)+(z) $ is an Artinian reduction of $X. $ Notice that $h_t= \Delta \HF_X(t)= \HF_X(t) -\HF_X(t-1) $ is the Hilbert function in degree $t$ of any Artinian reduction of $X $   and $h_t =0 $ for every $t> s $ where $s $ is called  the {\bf{ socle degree}} of $R/I(X).$  Since $R/I(X)$ is Cohen-Macaulay, we recall that  $s$ coincides with the  regularity of $R/I(X).$   The vector $(h_0, \dots, h_s) $ is the $h$-vector of $R/I(X) $ or, for short, the $h$-vector of $X.$

Given a point $P=(a_0,\dots , a_n)\in \mathbb P^n_{\res}$ we define the dual linear form $L:=a_0 y_0+\cdots +a_n y_n$ in $\Gamma $ and we say that $L $ is the linear form associated to $P,$ actually it generates  the dual of $I(P).$

We present   a well known  result  concerning the inverse system of the ideal of a set of points.

\medskip
\begin{proposition}[\cite{IK99}, Theorem 1.15(2)]
\label{fat-Geramita}
Let $X=\{P_1,\dots, P_r\}$ be a set of distinct points of $\mathbb P^n_{\res}$ and let $L_1, \dots, L_r$ be the associated linear forms.
Then for all $j\ge 0$
$$
(I(X)^\perp)_j=
\langle L_1^{j},\dots ,L_r^{j}\rangle_{\res} .
$$
\end{proposition}
\medskip
Comparing Proposition \ref{fat-Geramita} and Proposition \ref{radical}, it is clear that in the case of zero dimensional schemes $X=\{P_1, \dots, P_r\},$ the a.f.g. $R$-submodules of $I^{\perp}$   are the inverse systems of $I(P_i),$ that is $M_i=\langle L_i^{j} : j\ge 0 \rangle.  $

\medskip

 By using  Theorem 3.14  and Remark 3.15 in \cite{ER22},    we characterize zero-dimensional Gorenstein schemes which are reduced.
The new part in the following result is $(3)$ implies $(1).$

\medskip
\begin{theorem}
\label{ger-link}
Let $X\subset \mathbb P^n_{\res}$ be  zero-dimensional Gorenstein scheme of degree $r$.
We denote by $s$ the socle degree of $R/I(X)$.
Let $I(X)^\perp=\langle H_t;t\ge 1\rangle$ be a $G$-admissible system of generators with respect to a linear form $z \in R$.
The following conditions are equivalent:
\begin{enumerate}
\item $X$ is reduced,
\item there  exist  (unique up scalars) linear forms $L_1,\dots, L_r\in \Gamma$  pairwise linear independent and  unique elements
$\alpha_1,\dots, \alpha_r\in \res$  such that for all $t\ge r+1$
    $$
    H_t=\frac{1}{(t+s-1)!}\sum_{i=1}^r \frac{\alpha_i}{z(P_i)} L_i^{t+s-1},
    $$
\item there  exist (unique up to scalars) linear forms $L_1,\dots, L_r\in \Gamma $ pairwise linear independent and  unique  elements
$\alpha_1,\dots, \alpha_r\in \res$ such that
    $$
    H_{r+2}=\frac{1}{(r+s+1)!}\sum_{i=1}^r \frac{\alpha_i}{z(P_i)} L_i^{r+s+1}.
    $$
\end{enumerate}
Under these conditions, $X=\{P_1,\dots,P_r\}$ where $P_i$ is the point defined by the linear form $L_i$, $i=1,\dots,r$, and  $\alpha_1,\dots, \alpha_r\in \res^*$.
\end{theorem}
\begin{proof}
Recall that $s=\deg(H_1)$, $\deg(H_t)=s+t-1$ and $r=\dim_{\res}\langle H_1\rangle$.

Now $(1)$ implies $(2)$ (and hence trivially $(3)$) follows by Theorem 3.13 and Remark 3.14 in \cite{ER22}.
Assume now $(3)$ and we prove that it implies $(1).$
Let $P_i$ be the point of $\mathbb P_{\res}^n$ defined by the linear form  $L_i$, $i=1,\dots,r$.
Since the linear forms are pairwise linear independent we have that $P_i\neq P_j$ if $i\neq j$.
We consider the set of points $Y=\{P_1,\dots, P_r\}$,
we will prove that $X=Y$.

From $(3)$ and \propref{fat-Geramita} we get
$$
\langle H_{r+2} \rangle \subset \langle L_1^{r+s-1},\dots L_r^{r+s-1}\rangle \subset I(Y)^\perp,
$$
so $I(Y)\subset \ann(H_{r+2}).$
Furthermore, by \cite[Proposition 2.11]{ER21} we deduce
$$
I(Y)_{\le r+1} R\subset \ann(H_{r+2})_{\le r+1} R=I(X).
$$

On the other hand, since the Castelnuovo-Mumford regularity of $I(Y)$ is at most $r$ we get
that $I(Y)_{\le r+1} R=I(Y)$.
Hence we have that
$$
I(Y)\subset I(X).
$$

Since the ground field $\res$ is infinite there is a linear regular element $x \in R_1$  of $R/I(X)$ and $R/I(Y)$.
We have $I(Y)+(x)\subset I(X)+(x)$.
Recall that  both $R/I(X)$ and $R/I(Y)$ are Cohen-Macaulay rings, so
$$
\length(R/I(X)+(x))=\length{(R/I(Y)+(x)})=r
$$
and then $I(Y)+(x)= I(X)+(x)$.
In particular the $h-$vector of $R/I(X)$ and $R/I(Y)$
coincides.
Being both rings graded and Cohen-Macaulay we get
that
$$
\HF_{R/I(Y)}=\HF_{R/I(X)}.
$$
Since $I(Y)\subset I(X)$ we deduce that $Y=X$.

Finally, if one the equivalent conditions holds, then $\alpha_1, \dots, \alpha_r \in \res^*$ by Theorem 3.13 in \cite{ER22}.
\end{proof}

 \bigskip
The condition $(3)$ in \thmref{ger-link} translates the problem to be  reduced for $X$ in terms of the identifiability of a specific form, see \cite{AC20}.
This means that we have to prove that $H_{r+2} $ can be uniquely written  as sum of powers of $r$ linear forms $L_1, \dots,L_r. $    If we know $\deg{X}=r $ and the socle degree $s$  of $R/I(X),$   then to verify  if $X$ is reduced it consists in an effective   computation on $\langle H_{r+2}\rangle  = (I(X) +(z^{r+2}))^{\perp} $, where $z$ is a linear regular element of $R/I(X)$.

 \medskip
In the following example we consider a reduced zero-dimensional Gorenstein scheme  $X\subset \mathbb P^2_{\res}$      of degree $r,$ hence the corresponding $H_{r+2}$ is given by  Theorem \ref{ger-link}(3).

 \vskip 2mm
\begin{example}
Let us consider the ideal $I$ of $R=\res[x_1,x_2,x_3]$ generated by the forms
$x_1^2-x_1x_3, x_2^2-x_2x_3$.
The projective variety $X=V(I)$ is a complete intersection of $\mathbb P_{\res}^2$ defining four points: $(1,0,1)$, $(0,1,1)$, $(0,0,1)$ and $(1,1,1)$.
Hence $X$ is a reduced Gorenstein (in fact a complete intersection)  set of $r=4$ points.
We can deduce this fact from the last result starting from the defining ideal $I$.
It is easy   to prove that $x_3$ is a non-zero divisor of $A=R/I, $ the socle degree $s=2, $ $r+2=6$ and that $H_{6}$, i.e. a generator of $(I+(x_3^6))^{\perp}$, can be computed by \cite{E-InvSyst14} :
\begin{multline*}
H_6=y_1^6y_2+3y_1^5y_2^2+5y_1^4y_2^3+5y_1^3y_2^4+3y_1^2y_2^5+y_1y_2^6+6y_1^5y_2y_3+15y_1^4y_2^2y_3+20y_1^3y_2^3y_3+\\15y_1^2y_2^4y_3+6y_1y_2^5y_3+15y_1^4y_2y_3^2+30y_1^3y_2^2y_3^2+
30y_1^2y_2^3y_3^2+15y_1y_2^4y_3^2+20y_1^3y_2y_3^3+\\30y_1^2y_2^2y_3^3+20y_1y_2^3y_3^3+15y_1^2y_2y_3^4+15y_1y_2^2y_3^4+6y_1y_2y_3^5
\end{multline*}
If $L_1=y_1+y_3$, $L_2=y_2+y_3$, $L_3=y_3$ and $L_4=y_1+y_2+y_3$ then
$$
H_6=\frac{1}{7}(-L_1^7-L_2^7+L_3^7+L_4^7)
$$
so condition $(3)$ of \thmref{ger-link} holds true.
Indeed  $X$ is a reduced Gorenstein zero-dimensional scheme.
\end{example}

\medskip

In the following example we consider a non-reduced zero-dimensional Gorenstein scheme  $X\subset \mathbb P^2_{\res}$      of degree $r,$ hence the corresponding $H_{r+2}$ cannot verify condition (3) of Theorem \ref{ger-link}.

\vskip 2mm
\begin{example}
Let us consider the ideal $I$ of $R=\res[x_1,x_2,x_3]$ generated by the forms $x_1^2+x_2^2-x_3^2, x_1^2-x_2 x_3-x_3^2$.
The projective variety $X=V(I)$ is a complete intersection of $\mathbb P_{\res}^2$ defining four points: a double point $(0,-1,1)$ and two simple points: $(1,0,1)$ and $(-1,0,1)$.
Hence $X$ is a non-reduced Gorenstein set of points.
In this case $r=4$ and $s=2$.
We have that $x_3$ is a non-zero divisor of $A=R/I$ and that $H_6$, i.e. a generator of $(I+(x_3^6))^{\perp}$ can be computed by \cite{E-InvSyst14}:
\begin{multline*}
    H_6=y_1^7-7 y_1 y_2^6+42 y_1 y_2^5 y_3+21 y_1^5 y_3^2-105 y_1 y_2^4 y_3^2\\
    +140 y_1 y_2^3 y_3^3+35 y_1^3 y_3^4-105 y_1 y_2^2 y_3^4+42 y_1 y_2 x_3^5
\end{multline*}

\noindent
A linear algebra  computation by \cite{DGPS}  shows that there are not $\alpha_1,\dots ,\alpha_4\in \res$
and  linear forms $L_1,\dots , L_4$ such that
 $$H_{6}=\frac{1}{7!}\sum_{i=1}^4 \alpha_i L_i^{7}.$$
 Hence we recover that $I$ is not radical from the condition $(3)$ of \thmref{ger-link}.
\end{example}


\providecommand{\bysame}{\leavevmode\hbox to3em{\hrulefill}\thinspace}
\providecommand{\MR}{\relax\ifhmode\unskip\space\fi MR }
\providecommand{\MRhref}[2]{%
  \href{http://www.ams.org/mathscinet-getitem?mr=#1}{#2}
}
\providecommand{\href}[2]{#2}

\end{document}